\documentclass[matann,numbook,envcountsame]{svjour}

\usepackage{amsmath,latexsym,amssymb,times,mathptm}

\def\sqr#1#2{{\vcenter{\hrule height.#2pt
        \hbox{\vrule width.#2pt height#1pt \kern#1pt
                \vrule width.#2pt}
        \hrule height.#2pt}}}
\def\square{\mathchoice\sqr64\sqr64\sqr{4}3\sqr{3}3}
\def\QED{\hfill$\square$}

\def\tratto{\mbox{\rule{2mm}{.2mm}$\;\!$}}

\begin{document}

\title{A formula for the core of an ideal\thanks{The
authors were partially supported by the NSF.}}

\titlerunning{A formula for the core of an ideal}
\authorrunning{C. Polini and B. Ulrich}

\author{Claudia Polini \and Bernd Ulrich}

%\offprints{}

\institute{
{\sc C. Polini} \\
Department of Mathematics, University of Notre Dame, Notre Dame,
IN 46556, USA \\ (e-mail$\colon$ cpolini@nd.edu) \and
{\sc B. Ulrich} \\
Department of Mathematics, Purdue University, West Lafayette,
IN 47907, USA \\
(e-mail$\colon$ ulrich@math.purdue.edu)}

\date{Received:  / Revised version: }

\maketitle

\begin{abstract}
The core of an ideal is the intersection of all its reductions.
For large classes of ideals $I$ we explicitly describe the core as
a colon ideal of a power of a single reduction and a power of $I$.

\medskip

\noindent \subclassname{Primary 13B22; Secondary 13A30, 13B21,
13C40, 13H10.}
\end{abstract}

\section{Introduction}

The purpose of this paper is to prove a formula for the core of an
ideal that had been conjectured by the authors and A. Corso in
\cite{CPU2}. For ideals $J \subset I$ in a Noetherian ring one
says that $J$ is a {\it reduction} of $I$, or $I$ is {\it
integral} over $J$, if $I^{n+1}=JI^n$ for some $n \geq 0$. The
{\it core} of $I$, ${\rm core}(I)$, defined as the intersection of
all reductions of $I$, is a somewhat mysterious subideal of $I$
that encodes information about the possible reductions of $I$. The
concept was introduced by Rees and Sally (\cite{RS}), and has been
studied further by Huneke and Swanson and by Corso and the authors
(\cite{HS}, \cite{CPU}, \cite{CPU2}, \cite{CPU3}). It has a close
relation to Brian\c{c}on-Skoda type theorems and to coefficient,
adjoint and multiplier ideals (\cite{HS}, \cite{Hy}, \cite{HyS}).
Moreover, Hyry and Smith recently discovered an unexpected
connection with Kawamata's conjecture on the non-vanishing of
sections of line bundles.
%by proving that Kawamata's conjecture would follow from a suitable
%graded analogue of the above conjecture on the core. They also
%prove a version of our original conjecture under assumptions that
%are natural in the geometric context (\cite{HyS}).
They showed that Kawamata's conjecture would follow from a formula
that essentially amounts to a graded analogue of the above
conjecture on the core. They were also able to prove our original
conjecture under additional assumptions that arise naturally in
the geometric context (\cite{HyS}).

If $R$ is a Noetherian local ring with infinite residue field $k$,
then ${\rm core}(I)$ is the intersection of {\it minimal}
reductions of $I$, i.e., of reductions minimal with respect to
inclusion. The minimal number of generators of every minimal
reduction of $I$ is the {\it analytic spread} of $I$, which can
also be defined as $\ell(I)= {\rm dim} \, {\rm gr}_I(R) \otimes_R k$.
Given a reduction $J$ of $I$ we write $r_J(I)$ for the least
integer $n \geq 0$ such that $I^{n+1} = J I^n$, and we define the
{\it reduction number} $r(I)$ of $I$ to be ${\rm min}\{ r_J(I)\}$,
where $J$ ranges over all minimal reductions of $I$. In this paper
we prove the following result:

Let $(R, {\mathfrak m})$ be a local Gorenstein ring with infinite
residue field $k$, let $I$ be an $R$-ideal with $g={\rm ht} \, I
> 0$ and $\ell= \ell(I)$, and let $J$ be a minimal reduction of $I$
with $r=r_J(I)$. Assume $I$ satisfies $G_{\ell}$, ${\rm depth} \,
R/I^j \geq {\rm dim} \, R/I -j+1$ for $1 \leq j \leq \ell-g$, and
either ${\rm char} \,  k =0 $ or ${\rm char} \, k  > r -\ell + g
$. Then
\[
{\rm core}(I)= J^{n+1} : I^n
\]
for every $n \geq {\rm max} \{ r - \ell + g, 0 \}$.

We are now going to discuss the assumptions in this theorem. The
condition on the characteristic is vacuous if $r \leq \ell -g +1$,
in which case the result has been shown in \cite{CPU2}. The
$G_{\ell}$ property is a rather weak requirement on the local
number of generators of $I$, that is always satisfied if
$I_{\mathfrak p}$ can be generated by ${\rm dim} \, R_{\mathfrak
p}$ elements for every prime ideal ${\mathfrak p} \not= {\mathfrak
m}$ containing $I$. Both assumptions, the $G_{\ell}$ condition and
the inequalities ${\rm depth} \,  R/I^j \geq {\rm dim} \, R/I -j
+1 $ for $1 \leq j \leq \ell-g$, are automatically satisfied if
$I$ is ${\mathfrak m}$-primary, or more generally, if $I$ is {\it
equimultiple}, i.e., $\ell = g$. They also hold for
one-dimensional generic complete intersection ideals, or more
generally, for Cohen-Macaulay generic complete intersections with
{\it analytic deviation one}, i.e., $\ell=g+1$. In the presence of
the $G_{\ell}$ property, the depth inequalities for the powers
obtain if $I$ is perfect with $g=2$, $I$ is perfect Gorenstein
with $g=3$, or more generally, if $I$ is in the linkage class of a
complete intersection.

Although most of the paper is devoted to proving the main theorem,
 Theorem~\ref{main}, we show some of the auxiliary results in greater generality than
needed since they might be interesting in their own right. For
instance, to prove the inclusions $J^{n+1} : I^n \subset {\rm
core}(I)$, we identify the first ideals with graded components of
the canonical module of the {\it extended Rees algebra} $R[It,
t^{-1}]\subset R[t, t^{-1}]$ of $I$. The computation of the
canonical module is the content of Proposition~\ref{canonical}. It
relies on residual intersection theory and earlier results from
\cite{U2}, where graded components of canonical modules have been
used to identify colon ideals of powers of reductions that are
independent of the chosen reduction. The computation of graded
components of canonical modules has also played an important role
in \cite{Hy}, \cite{CPU2} and \cite{HyS}. To prove the reverse
containment ${\rm core}(I) \subset J^{n+1} : I^n$ we first
consider the case $\ell=g=1$. In this case Theorem~\ref{heightone}
provides a general formula for ${\rm core}(I)$ that holds in a
Cohen-Macaulay ring of arbitrary characteristic and specializes to
the desired equality ${\rm core}(I)=J^{n+1} : I^n$ under suitable
assumptions on the characteristic. Important ingredients in this
proof are ideas from \cite{HS} as well as a use of characteristic
zero that was inspired to us by K. Smith. In Lemma~\ref{lifting}
and Theorem~\ref{submain} we lift the containment obtained for
$\ell=g=1$ to establish an inclusion for the core of ideals of
arbitrary height and arbitrary analytic spread. Once
Lemma~\ref{lifting}, the main new technical result, is in place,
the proof of Theorem~\ref{submain} essentially follows the lines
of \cite{CPU2}. For both results we need to use residual
intersection theory since we are not restricting ourselves to the
case of equimultiple ideals. We give various classes of examples
showing that the main theorem may fail if \, $0 < {\rm char} \, k
\leq r - \ell +g$ or if any of the other assumptions are dropped.

The formula ${\rm core}(I)= J^{n+1} : I^n $ has been proved
independently by Huneke and Trung under the assumption that $R$ is
a local Cohen-Macaulay ring whose residue field has characteristic
zero and $I$ is an equimultiple ideal (\cite{HT}). Both papers,
the present one and \cite{HT}, were preceded by the work of Hyry
and Smith who proved the same formula if in addition the Rees ring
of $I$ is Cohen-Macaulay (\cite{HyS}). We would like to thank
Karen Smith for sharing her ideas with us before \cite{HyS} was
completed.

\section*{Preliminaries}

We begin by reviewing some definitions and basic facts. Let $R$ be
a Noetherian ring, $I$ an $R$-ideal and $j,s$ integers. We set
$I^j=R$ whenever $j \leq 0$. The ideal $I$ satisfies condition
$G_s$ if for every prime ideal ${\mathfrak p}$ containing $I$ with
${\rm dim} \, R_{\mathfrak p} \leq s-1$, the minimal number of
generators $\mu(I_{\mathfrak p})$ of $I_{\mathfrak p}$ is at most
${\dim} \, R_{\mathfrak p}$. One says that $I$ satisfies
$G_\infty$ in case $G_s$ holds for every $s$. When writing that
${\mathfrak a} \, \colon I$ is a {\it geometric $s$-residual
intersection} of $I$ we mean that ${\mathfrak a}$ is an
$s$-generated $R$-ideal properly contained in $I$ and that ${\rm
ht} \, {\mathfrak a}  \, \colon I \geq s $, ${\rm ht} (I,
{\mathfrak a} \, \colon I) \geq s+1$. The ideal $I$ is said to be
of {\it linear type} if the natural map from the symmetric algebra
to the Rees algebra of $I$ is an isomorphism; in this case $I$ has
no proper reductions.

Now assume in addition that $R$ is a local Cohen-Macaulay ring,
and let $H_{\bullet}$ denote the homology of the Koszul complex on
a generating sequence $f_1, \ldots ,f_n$ of $I$. One says that $I$
satisfies {\it sliding depth} if ${\rm depth} \, H_i \geq {\rm
dim} \, R -n +i$ for every $i \geq 0$, and that $I$ is {\it
strongly Cohen-Macaulay} if $H_{\bullet}$ is Cohen-Macaulay. These
notions fit into the following sequence of implications: In case
$I$ is a perfect ideal of grade $2$ or a perfect Gorenstein ideal
of grade $3$, then $I$ is in the linkage class of a complete
intersection, which in turn implies that $I$ is strongly
Cohen-Macaulay (\cite[1.11]{H1}). Strong Cohen-Macaulayness
obviously implies the sliding depth property. If $I$ satisfies
sliding depth and $G_\infty$, then $I$ is of linear type and the
associated graded ring ${\rm gr}_I(R)$ is Cohen-Macaulay
(\cite[9.1]{HSV}). If on the other hand $I$ is strongly
Cohen-Macaulay of ${\rm height} \, g$ and satisfies $G_s$, then
${\rm depth} \,  R/I^j \geq {\rm dim} \, R/I -j +1 $ for $1 \leq j
\leq s-g+1$ (\cite[the proof of 5.1]{HSV}). The latter condition
in turn implies that $R/ {\mathfrak a}\, \colon I$ is
Cohen-Macaulay for every (geometric) $s$-residual intersection of
$I$, at least if $R$ is Gorenstein and $I$ satisfies $G_s$
(\cite[2.9$({\it a})$]{U}).

\section{The canonical module}

We begin by computing the graded canonical module of certain
extended Rees rings.

\begin{proposition}\label{canonical}
Let $R$ be a Noetherian local ring, $I$ an $R$-ideal with ${\rm
ht} \, I > 0$, and $J$ a reduction of \, $I$ with $r= r_J(I)$.
Write $A=R[Jt,t^{-1}] \subset B=R[It,t^{-1}]$. Assume that $A$ is
Cohen-Macaulay with graded canonical module $\omega_A \cong (1,
It)^s A(a)$ for some integers $s$ and $a$. Then for every integer
$n \geq {\rm max} \{ r-s, 0 \}$,
\[
\omega_B \cong (A :_{R[t,t^{-1}]} I^{n}) (a+n). \]
\end{proposition}
\begin{proof}
First notice that $R$ is Cohen-Macaulay. We write $K= {\rm
Quot}(R)$ and make the identification $\omega_A = (1, It)^s
At^{-a} \subset R[t,t^{-1}]$. For $L=(1,(1, It)^sI^{n}t^{n}) A
\subset R[t,t^{-1}]$ consider the exact sequence of $A$-modules,
\[
0 \rightarrow L \longrightarrow B \longrightarrow C \rightarrow 0.
\]
\noindent Since $s+n \geq r_J(I)$ it follows that $C$ is
concentrated in finitely many degrees, and hence has grade $\geq
2$. Thus, dualizing the above exact sequence into $\omega_A$ we
obtain
\begin{eqnarray*}
\omega_B  & \cong & {\rm Hom}_A(B, \omega_A) \cong {\rm Hom}_A(L, \omega_A) \\
& \cong & \omega_A :_{K(t)} L\\
& = & \omega_A \cap (\omega_A :_{K(t)} (1,It)^sI^{n}t^{n}A)\\
& = & \omega_A \cap (\omega_A :_{K(t)} \omega_A I^{n}t^{a+n})\\
& = & \omega_A \cap ((\omega_A :_{K(t)} \omega_A) :_{K(t)} I^{n}t^{a+n})\\
& = & \omega_A \cap (A :_{K(t)} I^{n}t^{a+n})\\
& = & A :_{\omega_A} I^{n}t^{a+n}.
\end{eqnarray*}
Finally, for every integer $i$ one has
\begin{eqnarray*}
[A :_{R[t,t^{-1}]} I^{n}t^{a+n}]_i  & = & (J^{i+a+n} :_R I^{n})t^i \\
& \subset & (J^{i+a+n} :_R J^{n})t^i\\
& = & J^{i+a}t^i,
\end{eqnarray*}
\noindent where the last equality holds because ${\rm gr}_J(R)$ is
Cohen-Macaulay, ${\rm ht} \, J > 0$, and $n \geq 0$. However,
\[
J^{i+a}t^i = [At^{-a}]_i \subset [\omega_A]_i.
\]
\noindent Therefore $A :_{R[t,t^{-1}]} I^{n}t^{a+n} \subset
\omega_A \subset R[t,t^{-1}],$ or equivalently,
\[
A :_{\omega_A} I^{n}t^{a+n}= A :_{R[t,t^{-1}]} I^{n}t^{a+n}.
\]
\noindent Thus by the above,
\[
\omega_B \cong (A :_{R[t,t^{-1}]} I^{n})(a+n)
\]
as claimed. \QED
\end{proof}

\medskip

\begin{remark}\label{remarkcanonical}{\rm

\

\begin{enumerate}
\item
In Proposition~\ref{canonical}, the condition that $A$ is
Cohen-Macaulay can be replaced by the weaker assumption that $A$
satisfies $S_2$ and $R$ is universally catenary.

\item Let $R$ be a local Gorenstein ring with infinite residue
field, let $I$ be an $R$-ideal with $g={\rm ht} \, I > 0$ and
$\ell= \ell(I)$, and assume that $I$ satisfies $G_{\ell}$ and
${\rm depth} \, R/I^j  \geq {\rm dim} \, R/I -j+1$ for $1 \leq j
\leq \ell-g$. According to \cite[the proof of 2.1]{U2}, every
minimal reduction $J$ of $I$ satisfies the assumptions of
Proposition~\ref{canonical} with $s= \ell -g$ and $a=1-g$. Thus
for every $n \geq {\rm max} \{r_J(I)-\ell + g, 0 \}$,
\[
\omega_B=(A :_{R[t, t^{-1}]} I^{n})t^{g -n-1}
\]
is a graded canonical module of $B$.

\noindent In particular, for every $i$ and every $n \geq {\rm max}
\{ r_J(I) -\ell +g, 0\}$,
\[
[\omega_B]_{i+g-1}t^{-i-g+1} = J^{i+n} :_R I^{n}.
\]
\end{enumerate} }
\end{remark}

\bigskip

\begin{corollary}\label{balanced}
Let $R$ be a local Gorenstein ring with infinite residue field,
let $I$ be an $R$-ideal with $g={\rm ht} \, I > 0$ and $\ell=
\ell(I)$, and assume that $I$ satisfies $G_{\ell}$ and ${\rm
depth} \, R/I^j  \geq {\rm dim} \, R/I -j+1$ for $1 \leq j \leq
\ell-g$. Then for every fixed integer $i$, the ideal
\[
 J^{i+n} :_R I^{n}
\]
is independent of $J$ and $n$, as long as $J$ is a minimal
reduction of $I$ and $n \geq {\rm max} \{r_J(I) - \ell +g, 0\}$.
\end{corollary}
\begin{proof}
By Remark~\ref{remarkcanonical}(2) one has $\omega_B=(A :_{R[t,
t^{-1}]} I^{n})t^{g-n-1}$ and it suffices to show that this
submodule of $R[t, t^{-1}]$ is uniquely determined by $B$. Notice
that $\omega_B$ is a graded canonical module of $B$ and a graded
submodule of $R[t, t^{-1}]$, and that $[\omega_B]_i=Rt^i$ for $i
\ll 0$. Now the first two properties determine this submodule up
to multiplication with a unit $u$ in ${\rm Quot}(R)$, and the last
property then forces $u$ to be a unit in $R$. \QED
\end{proof}

\smallskip

The next result has been shown in \cite[3.2 and 3.4]{Hy} under the
assumption that $R[It]$ is Cohen-Macaulay.

\smallskip

\begin{remark}\label{Jout}{\rm
In addition to the assumptions of Corollary~\ref{balanced} assume
that ${\rm gr}_I(R)$ is Cohen-Macaulay. Then for every $i \geq 0$ and
every $n \geq {\rm max} \{r_J(I) - \ell +g, 0\}$,
\[
 J^{i+n} : I^n =J^i (J^n :I^n) = I^i (J^n : I^n).
\]
In particular $J^n \colon I^n$ is the coefficient ideal of $I$
with respect to $J$ in the sense of \cite[2.1]{AH2}.

Indeed, the claim about the coefficient ideal ${\mathfrak a}$
follows from the second asserted equality because it gives $J^n \colon
I^n \subset {\mathfrak a}$, whereas the reverse inclusion is
always true. To prove the equalities we may assume that $R$ is
complete. By Remark~\ref{remarkcanonical}(2) it suffices to
show that $\omega_B$ as a graded $A$-module
is generated in degrees $\leq g-1$.
%After factoring out part of a general system of
%parameters of $R$ we may assume that $\ell=d={\rm dim}(R)$. Now
%${\rm pd}_R(I^j/I^{j+1}) \leq j+g$ for every $j \geq 0$. on the
%other hand,
Let $T$ be a regular local ring mapping onto $R$, set $c={\rm dim}
\, T - {\rm dim} \, R$, and consider the polynomial ring $S=T[X_1,
\ldots, X_{\ell}]$. Mapping $S$ homogeneously onto ${\rm
gr}_J(R)$, the ring ${\rm gr}_I(R)$ becomes a finite graded
$S$-module whose homogeneous minimal free resolution $F_{\bullet}$
has length $c + \ell $. As ${\rm pd}_T(I^j/I^{j+1}) \leq c+g+j <
c+ \ell$ for every $ 0 \leq j \leq \ell -g -1$, it follows that
$F_{c+\ell}$ is generated in degrees $\geq \ell-g$. Therefore the
canonical module of ${\rm gr}_I(R)$ is generated in degrees $\leq
g$ as a graded module over $S$, hence over ${\rm gr}_J(R)$. Thus
indeed $\omega_B$ is generated in degrees $\leq g-1$ as a graded
$A$-module. }
\end{remark}

\bigskip

\section{The case of analytic spread one}

The next lemma is a minor modification of \cite[2.2]{CPU2}, which
in turn was based on \cite[the proof of 3.8]{HS}. To simplify
notation we write $x \colon y$ instead of $(x) :_{R} (y)$ for
elements $x, y$ of a ring $R$.

\begin{lemma}\label{oldlemma}
Let $(R, {\mathfrak m}, k)$ be a Noetherian local ring, let $K$ be
an $R$-ideal, and let $x$, $y$ be elements of $R$ such that $yK
\subset xK$ and $x$ is a non-zerodivisor. Let $c
> {\rm dim}_k (K\colon {\mathfrak m} /K)$ and
let $u_1, \ldots, u_c$ be units in $R$ that are not all congruent
modulo ${\mathfrak m}$. Then for every $j \geq 0$,
\[
x(K \colon {\mathfrak m}) \cap  \bigcap_{i=1}^{c}(x+u_i y)(K
\colon {\mathfrak m}) \subset x(x^j \colon y^j) \cap
\bigcap_{i=1}^{c} (x+u_i y)(x^j \colon y^j).
\]
\end{lemma}
\begin{proof}
Let $\alpha$ be an element of the intersection on the left
hand side. Write
\[
\alpha = xs= (x+u_1 y)s_1=  \ldots = (x+u_c y)s_c
\]
where $s$ and all $s_i$ belong to $K \colon {\mathfrak m}$. We are
going to prove by induction on $j$ that $s$ and all $s_i$ are in
$x^j \colon y^j$.

The assertion being trivial for $j=0$ we may assume that $j>0$. We
first show that $s_i \in x^j \colon y^j$ for $1 \leq i \leq c$. By
our induction hypothesis,
\[
s_i(x+u_iy)= \alpha \in x(x^{j-1} \colon y^{j-1}).
\]
Hence
\begin{eqnarray*}
s_i& \in & (x(x^{j-1} \colon y^{j-1})) \colon (x+u_iy)  \\
& \subset &  (x^{j} \colon y^{j-1}) \colon (x+u_iy) \\
& = &  x^{j} \colon (x+u_iy)y^{j-1}  \\
& = &  x^{j} \colon (xy^{j-1}+u_iy^j).
\end{eqnarray*}
Since $s_i \in x^{j-1} \colon y^{j-1}=  x^j \colon xy^{j-1}$ by
induction hypothesis, it follows that $s_i \in x^j \colon u_iy^j$.
Therefore $s_i \in x^j \colon y^j$, as asserted.

Next we prove that $s \in x^j \colon y^j$. As $yK \subset xK$ one
has $K \subset x^j \colon y^j$. If $s_i \in K$ for some $i$, then
\[
xs=(x+u_iy)s_i \in xK +yK=xK \subset x(x^j \colon y^j).
\]
Thus $s \in x^j \colon y^j$ since $x$ is a non-zerodivisor. So we
may assume that $s_i \not \in K$ for $1 \leq i \leq c$. Let
`${}^{\tratto}$' denote images in $\overline{R}=R/K$. Now
$\overline{s}_1, \ldots, \overline{s}_c$, or equivalently,
$\overline{u}_1\overline{s}_1, \ldots,
\overline{u}_c\overline{s}_c$ are $c$ nonzero elements of the
$k$-vector space $\overline{K \colon {\mathfrak m}}$. They are
linearly dependent over $k$, and after shrinking $c$ if needed we
may assume that $c$ is minimal with respect to this property and
that $\overline{u}_1, \ldots, \overline{u}_c$ are still not all
equal. Obviously $c \geq 2$. Now there exist units $\lambda_1,
\ldots, \lambda_c$ in $R$ so that $\displaystyle\sum_{i=1}^c
\overline{\lambda}_i \overline{u}_i \overline{s}_i =
\overline{0}$. Notice that $\displaystyle\sum_{i=1}^c
\overline{\lambda}_i \overline{s}_i \not= \overline{0}$ because
$\overline{u}_1, \ldots, \overline{u}_c$ are not all equal. Since
$\displaystyle\sum_{i=1}^c \lambda_i u_i s_i \in K$, the element
$\displaystyle\sum_{i=1}^c y \lambda_i u_i s_i $ belongs to $yK
\subset xK$, hence can be written as $x \xi $ for some $\xi \in
K$. Set $\lambda = \displaystyle \sum_{i=1}^c \lambda_i$ and
multiply both sides by $\alpha$. Rewriting $\alpha$ by means of
the above equations, cancelling $x$, and using the containments
$s_i \in x^j \colon y^j$, we obtain
\[
\lambda s = \sum_{i=1}^c \lambda_i s_i + \xi \in (x^j \colon y^j)
+ K = x^j \colon y^j.
\]
If $\lambda \in {\mathfrak m}$ then $\lambda s \in K$ since $s \in
K \colon {\mathfrak m}$, and we conclude that $\overline{0} =
\displaystyle \sum_{i=1}^c \overline{\lambda}_i \overline{s}_i$,
which is impossible. Thus $\lambda$ is a unit, and the desired
inclusion $s \in x^j \colon y^j$ follows. \QED
\end{proof}

\medskip

\begin{lemma}\label{help}
Let $R$ be a local Cohen-Macaulay ring with infinite residue
field, let $I$ be an $R$-ideal with $\ell(I)= {\rm ht} \, I = 1$
and $r= r(I)$, let $J$ and $H$ be minimal reductions of $I$, and
let $n \geq r$ and $i$ be integers.
\begin{description}[$({\it a})$]
\item[$({\it a})$] $H^i(J^n : I^n)=I^i(J^n : I^n)$, and this ideal
is independent of $J$, $H$ and $n$. \item[$({\it b})$] $I(J^n :
I^n) \subset {\rm core}(I)$.
\end{description}
\end{lemma}
\begin{proof}
To prove $({\it a})$ write $K={\rm Quot}(R)$ and consider $S=
\bigcup_{j \geq 0} (I^j :_{K}I^j)$, the blowup ring of $I$. Let
$x$ be a generator of $J$ and notice that $r_J(I)=r$ (see, e.g.,
\cite[2.1]{Huc}). One has $S= I^n :_K I^n \subset I^n
:_KJ^n=I^n\frac{1}{x^n} \subset R[\frac{I}{x}] \subset S$, in
particular $I^n\frac{1}{x^n} = S$ \, (see also \cite[1.1 and its
proof]{Lip}). Hence $J^n:_RI^n = J^n :_K I^n = R:_K
I^n\frac{1}{x^n} =R:_K S$ is independent of $J$ and $n$.
Furthermore $IS$ is principal, hence $IS=HS$. As $J^n:_R I^n=R :_K
S$ is an $S$-ideal we conclude that
\begin{eqnarray*}
H^i(J^n :_R I^n) & = & H^iS(J^n :_R I^n)  \\
& = & I^iS(J^n :_R I^n)  \\
& = &  I^i(J^n :_R I^n).
\end{eqnarray*}
This proves $({\it a})$. Part $({\it b})$ follows because $({\it
a})$ implies that $I(J^n :_R I^n) \subset H$ for every choice of
$H$. \QED
\end{proof}

\medskip

We are now ready to prove our formula for the core of equimultiple
height one ideals. We first remark that in this case the core is the
intersection of a specific (finite) number of general principal
reductions. Here we denote the Hilbert-Samuel multiplicity by
$e(-)$ and length by $\lambda(-)$.

\smallskip

\begin{remark}\label{remarkblah}{\rm
Let $R$ be a local Cohen-Macaulay ring with infinite residue
field, let $I$ be an $R$-ideal with $\ell(I)= {\rm ht} \, I = 1$,
and write $t= \max \{ e(I_{\mathfrak p}) \, | \, {\mathfrak p} \in
{\rm Min}(I) \}$. Then ${\rm core}(I)$ is the intersection of $t$
general principal ideals in $I$.

Indeed,  \cite[4.9 and 4.5]{CPU} gives that ${\rm core}(I)$ is the
intersection of $ \max \{ {\rm type}(R_{\mathfrak p}/{\rm
core}(I_{\mathfrak p})) \, | \, {\mathfrak p} \in {\rm Min}(I) \}$
general principal ideals in $I$. On the other hand let ${\mathfrak
p} \in {\rm Min}(I)$ and let $J$ be a minimal reduction of
$I_{\mathfrak p}$. Write ${\rm core}(I_{\mathfrak p})= JL$ for
some $R_{\mathfrak p}$-ideal $L$. As ${\rm core}(I_{\mathfrak p})
: {\mathfrak p}_{\mathfrak p} = JL : {\mathfrak p}_{\mathfrak p}
\subset JL : J = L$, it follows that ${\rm type}(R_{\mathfrak
p}/{\rm core}(I_{\mathfrak p})) \leq \lambda(L/JL) =
\lambda(R_{\mathfrak p}/J)= e(I_{\mathfrak p})$.  }
\end{remark}

\smallskip

\begin{theorem}\label{heightone}
Let $R$ be a local Cohen-Macaulay ring with infinite residue field
$k$, let $I$ be an $R$-ideal with $\ell(I)= {\rm ht} \, I = 1$ and
$r= r(I)$, and let $J$ be a minimal reduction of $I$. Let $(y_1), \ldots,(y_t)$
be minimal reductions of $I$ so that ${\rm core} (I) = (y_1) \cap \cdots \cap
(y_t)$ and write $s= {\rm max} \{ r((J, y_i)) \, | \, 1 \leq i \leq
t\}$.
% Write $t=
%{\rm max} \, \{ e(I_{\mathfrak p}) \, | \, {\mathfrak p} \in {\rm
%Min} \, (I) \}$, let $ y_1, \ldots, y_t$ be general elements of
%$I$, and set $s= {\rm max} \{ r((J, y_i)) \, | \, 1 \leq i \leq
%t\}$.
\begin{description}[$({\it a})$]
\item[$({\it a})$]

\

\vspace{-.5cm}

\noindent $\begin{array}{ccc} {\rm core} (I)& = & J^{n+1} :
\displaystyle \sum_{y \in I} (J, y)^n
=  J ( J^n : \displaystyle \sum_{y \in I} (J, y)^n) \\
& = &  J^{n+1} : \displaystyle \sum_{i=1}^{t} (J, y_i)^n  =  J (
J^n : \displaystyle \sum_{i=1}^{t} (J, y_i)^n)
\end{array}$

\noindent for every $n\geq s$.

\

\item[$({\it b})$]
If ${\rm char} \, k =0$ or  ${\rm char} \, k > r$, then
\[
{\rm core}(I) =J^{n+1} :I^n = J (J^n:I^n)
\]
for every $n \geq r$.
\end{description}
\end{theorem}
\begin{proof}
First notice that the second and fourth equality in $({\it
a})$ and the second equality in $({\it b})$ are obvious because
$J$ is generated by a single regular element.

We now prove part $({\it a})$.
Since $(J,y)$ is a reduction of $I$ and $(y_i)$ is a reduction of
$(J, y_i)$, it follows that
\begin{eqnarray*} {\rm core} (I)& \subset & \bigcap_{y\in I} {\rm core}((J,
y))  \subset \bigcap_{i=1}^t {\rm core}((J,y_i)) \\
& \subset &  \bigcap_{i=1}^{t}(y_i)  =  {\rm core}(I).
\end{eqnarray*}
Therefore
\[
{\rm core}(I)=  \bigcap_{y\in I} {\rm core}((J, y))  =
\bigcap_{i=1}^t {\rm core}((J,y_i)).
\]

On the other hand,
\[
\begin{array}{ccccc}
J^{n+1} : \displaystyle \sum_{y \in I} (J, y)^n  & \subset &
J^{n+1} : \displaystyle\sum_{i=1}^{t} (J, y_i)^n & & \\
\parallel & & \parallel & & \\
\displaystyle \bigcap_{y \in I} (J^{n+1} : (J, y)^n) & \subset &
\displaystyle \bigcap_{i=1}^t (J^{n+1} :
(J, y_i)^n) & \subset & \displaystyle \bigcap_{i=1}^t {\rm core}((J,y_i)),
\end{array}
\]
where the last containment follows from Lemma~\ref{help}$({\it
b})$, applied to the ideals $(J,y_i)$. Thus it suffices to prove
that ${\rm core}((J,y))\subset J^{n+1} : (J,y)^n$.

To this end we may assume that $n \geq r((J,y))$, because $J^{n+1}
:(J,y)^n$ form a decreasing sequence of ideals. Furthermore
for every associated prime ${\mathfrak p}$ of $J$, ${\rm dim} \,
R_{\mathfrak p} =1$ and $({\rm core}((J,y)))_{\mathfrak p}= {\rm
core}((J,y)_{\mathfrak p})$ by \cite[4.8]{CPU}. Thus after
localizing at ${\mathfrak p}$ we may suppose that ${\rm dim} \,
R=1$. Write ${\mathfrak m}$ for the maximal ideal of $R$, $J=(x)$
and $K=J^n:(J,y)^n$. Lemma~\ref{help}$({\it a})$, applied to the
ideal $(J,y)$, shows that $yK\subset xK$. We use the notation of
Lemma~\ref{oldlemma}, assuming in addition that $z_i=x+u_iy$
generate minimal reductions of $(J,y)$. Now $z_iK=xK$ by
Lemma~\ref{help}$({\it a})$, and hence
\begin{eqnarray*}
(z_i) \cap (xK : {\mathfrak m})  & = & (z_i) \cap (z_iK : {\mathfrak m})
=  z_i((z_iK : {\mathfrak m}) :z_i)  \\
& = &  z_i(z_iK : z_i{\mathfrak m})   =  z_i(K : {\mathfrak m})  .
\end{eqnarray*}
Likewise
\begin{eqnarray*}
(x) \cap (xK : {\mathfrak m})  & = &  x(K : {\mathfrak m})  .
\end{eqnarray*}
Using these facts and Lemma~\ref{oldlemma} we deduce
\begin{eqnarray*}
(x) \cap (z_1) \cap \ldots \cap (z_c) \cap (xK : {\mathfrak m}) &
= & x(K : {\mathfrak m}) \cap z_1(K : {\mathfrak m})
\cap \ldots \cap z_c(K : {\mathfrak m})\\
& \subset & x(\bigcap_{j=1}^n (x^j :y^j)) = x(x^n :(x,y)^n)= xK.
\end{eqnarray*}
As $xK$ is an ${\mathfrak m}$-primary ideal it follows that $(x)
\cap (z_1) \cap \ldots \cap (z_c) \subset xK$, hence ${\rm
core}((J,y)) \subset(x) \cap (z_1) \cap \ldots \cap (z_c) \subset
J^{n+1} : (J,y)^n $. This completes the proof of $({\it a})$.

To prove part $({\it b})$ notice that our assumption on the
characteristic gives $I^r = \sum_{y \in I} (J,y)^r$. Since
$r=r_J(I)$ (see, e.g., \cite[2.1]{Huc}) we obtain $I^j = \sum_{y
\in I}(J,y)^j$ for $j \gg 0 $. Thus by part  $({\it a})$, ${\rm
core}(I)= J^{j+1} : I^j$. However $J^{j+1} : I^j =J^{n+1} : I^n$
for every $n \geq r$. \QED
\end{proof}

\bigskip

\section{The proof of the main Theorem}

The next lemma, though elementary, plays an important role in the
proof of Lemma~\ref{lifting}. Its use was inspired to us by K.
Smith.

\begin{lemma}\label{regular}
Let $R$ be a ring and let $x_1, \ldots,x_n$ be elements in $R$
such that $x_i, x_j$ form a regular sequence for all $1 \leq i < j
\leq n$. Then $(x_1) \cap \ldots \cap (x_n)=(x_1 \cdot \ldots
\cdot x_n)$.
\end{lemma}

\begin{lemma}\label{lifting}
Let $R$ be a local Cohen-Macaulay ring with infinite residue field
and assume that $R$ has a canonical module. Let $J$ be an
$R$-ideal with $\ell=\mu (J) > 0 $ satisfying $G_{\infty}$ and
sliding depth, and write
\[
{\mathcal A}={\mathcal A}(J) = \left\{ {\mathfrak a}  \, | \,
{\mathfrak a} \, \colon J \  \mbox{\rm is a geometric $(\ell
-1)$-residual intersection and} \  \mu (J/ {\mathfrak a}) = 1
\right\}.
\]
\noindent Let $t$ be a positive integer and let $H$ be an
$R$-ideal satisfying ${\rm ht}(J, J^t : H) \geq \ell.$ Then
\[
H \cap \bigcap_{{\mathfrak a} \in {\mathcal A}} ({\mathfrak a},
J^t) \subset J^t.
\]
\end{lemma}
\begin{proof} We induct on $\ell$. If ${\ell} = 1$ then
${\mathcal A}=\{0\}$ and the assertion is clear. Hence we may
assume ${\ell} \geq 2$. Let $ b \in H$ and suppose that $ b \in
J^{j-1} \setminus J^j$ for some $j$ with $1 \leq j \leq t$. We are
going to prove that there exists an ideal ${\mathfrak a} \in
{\mathcal A}$ with $b \not\in ({\mathfrak a},J^t)$. For this we
may assume that $b \in J$.

Let `${}^{\tratto}$' denote images in $\overline{R} = R/ 0 :
J^{\infty}$. Notice that $J \cap (0 : J^{\infty})=0$ since $J$
satisfies $G_1$. Thus the canonical epimorphism $R
\twoheadrightarrow \overline{R}$ induces isomorphisms $J^n \cong
\overline{J}^n$ for every $n \geq 1$. Therefore $\overline{b} \in
\overline{J}^{j-1} \setminus \overline{J}^j$,
$\mu(\overline{J})=\mu(J)=\ell$, and every ideal in ${\mathcal
A}(\overline{J}) $ is of the form $\overline{{\mathfrak a}}$ for
some ${\mathfrak a} \in {\mathcal A} (J)$. One trivially has that
${\rm ht} (\overline{J},  \overline{J}^t : \overline{H}) \geq
\ell$ and $\overline{J}$ satisfies $G_{\infty}$. Finally,
$\overline{R}$ is Cohen-Macaulay and $\overline{J}$ has the
sliding depth property by \cite[3.6]{HVV}. Thus we may replace $J
\subset R$ by $\overline{J} \subset \overline{R}$ to assume that
${\rm ht} \, J > 0.$

Write $G = {\rm gr}_J(R)$, $G_+ = \bigoplus_{i >0} G_i$, and
$b^*=b+J^j \in [G]_{j-1}.$ By \cite[9.1]{HSV}, $J$ is of linear
type and $G$ is Cohen-Macaulay. Also notice that $b^* \not= 0.$
Thus $b^* \not \in Q$ for some primary component $Q$ of $0$ in
$G$. Write $P= \sqrt{Q}$ and let ${\mathfrak p}$ be the preimage
of $P$ in $R$. We claim that
\begin{equation}\label{EQ}
{\rm ht} \, (G_+, Q)/Q \geq 2.
\end{equation}
By Cohen-Macaulayness $P$ is a minimal prime of $G$, therefore
$\ell(J_{\mathfrak p})= {\rm dim} \, R_{\mathfrak p}.$ If ${\rm
dim} \, R_{\mathfrak p} < \ell$ then $H_{\mathfrak p} \subset
J_{\mathfrak p}^t.$ Since $b \in H$, it would follow that
$\frac{b^*}{1}=\frac{0}{1}$ in $G_{\mathfrak p}$, hence $b^* \in
Q$, contradicting the choice of $Q$. Therefore ${\rm dim} \,
R_{\mathfrak p} \geq \ell, $ which gives $\ell(J_{\mathfrak p})=
{\rm dim} \, R_{\mathfrak p} \geq \ell \geq 2.$ On the other hand
as $G/P $ is a positively graded domain, $(G_+, P)/P$ is a prime
ideal and we obtain ${\rm ht} \, (G_+, P) / P = {\rm ht} \, (G_+,
P)_{\mathfrak p} / P_{\mathfrak p}. $ Since $J$ is of linear type,
$G \otimes_R R_{\mathfrak p}/{\mathfrak p}R_{\mathfrak p}$ is a
domain as well. Hence $P_{\mathfrak p} = {\mathfrak p}G_{\mathfrak
p}$ because $P_{\mathfrak p}$ is a minimal prime ideal. Now we
conclude that
\begin{eqnarray*}
{\rm ht}\, (G_+, Q)/Q & = & {\rm ht} \, (G_+, P)/P \\
& = & {\rm ht} \, (G_+, P)_{\mathfrak p} / P_{\mathfrak p}\\
& = & {\rm ht} \, (G_+, {\mathfrak p}G)_{\mathfrak p} / {\mathfrak
p}G_{\mathfrak p} \\
& = & \ell(J_{\mathfrak p}) \geq 2.
\end{eqnarray*}
This completes the proof of $(\ref{EQ})$.

Write $M$ for the homogeneous maximal ideal of $G$, $A=(G/Q)_M$,
$N = ((G_+, Q)/Q)_M,$ and $B={\rm End}_A (\omega_A)$ for the
$S_2$-ification of $A.$ Notice that $A \hookrightarrow B$ since
$A$ is unmixed. Furthermore ${\rm ht} \, N B = {\rm ht} \, N \geq
2, $ where the equality follows from \cite[3.5$({\it b})$]{HH} and
the inequality is implied by $(\ref{EQ})$. Therefore ${\rm grade}
\, N B \geq 2.$ Since $0 \not= b^*A \subset A \subset B $ and $N
\subset {\rm Rad} (B),$ by Krull's intersection theorem there
exists an integer $n$ so that $b^*A \not\subset (NB)^n.$ Let $x_1,
\ldots, x_n $ be $n$ general elements of $J$ and write
$x_i^*=x_i+J^2 \in [G]_1$. As grade $NB \geq 2, $
Lemma~\ref{regular} implies that
\begin{eqnarray*}
x_1^*A \cap \ldots \cap x_n^*A & \subset & x_1^*B \cap \ldots \cap
x_n^*B
\\
& = & x_1^* \cdot \ldots \cdot x_n^* B\\
& \subset & (NB)^n.
\end{eqnarray*}
Therefore $b^*A \not\subset x_i^* A$ for some $i$. Write $x=x_i$
and let `${}^{\tratto}$' denote images in $\overline{R} = R/ (x).$
Notice that $b^* \not\in x^*G.$ By the general choice of $x$ and
since ${\rm ht} \, J > 0$, it follows that $x$ is $R$-regular.
Thus the Cohen-Macaulayness of $G$ and the genericity of $x$ imply
${\rm gr}_{\overline{J}}(\overline{R}) = G/x^*G$. Therefore
$\overline{b} \not\in \overline{J}^j.$ Again by the general choice
of $x$, $\mu(\overline{J})=\ell -1$ and $\overline{J}$ satisfies
$G_{\infty}.$ Moreover $\overline{J}$ has the sliding depth
property according to \cite[3.5]{HVV}, ${\rm ht} (\overline{J},
\overline{J}^t : \overline{H}) \geq \ell -1$ and every ideal of
${\mathcal A}(\overline{J}) $ is of the form $\overline{{\mathfrak
a}}$ for some ${\mathfrak a} \in {\mathcal A} (J).$ Now the
induction hypothesis shows that $\overline{b} \not\in
(\overline{{\mathfrak a}},\overline{J}^t)$ for some ${\mathfrak a}
\in {\mathcal A} (J),$ hence $b \not\in ({\mathfrak a},J^t)$. \QED
\end{proof}

\medskip

\begin{remark}\label{lineartype}{\rm
Let $R$ be a local Gorenstein ring with infinite residue field and
let $I$ be an $R$-ideal with $g={\rm ht} \, I$ and $\ell=
\ell(I)$. Assume that $I$ satisfies $G_{\ell}$ and ${\rm depth} \,
R/I^j \geq {\rm dim} \, R/I -j+1$ for $1 \leq j \leq \ell-g$. Then
any minimal reduction $J$ of $I$ satisfies $G_{\infty}$ and
sliding depth, as required in Lemma~\ref{lifting}. Furthermore $J$
is of linear type, ${\rm gr}_J(R)$ is Cohen-Macaulay, and ${\rm
ht} \, J \colon I \geq \ell$.

Indeed by \cite[2.9({\it a}) and 1.11]{U}, ${\rm ht} \, J : I \geq
\ell$. Hence $J$ satisfies $G_{\infty}$. Now according to
\cite[2.9({\it a}), 1.12 and 1.8({\it c})]{U} $J$ has the sliding
depth property. Therefore $J$ is of linear type and ${\rm
gr}_J(R)$ is Cohen-Macaulay by \cite[9.1]{HSV}. }
\end{remark}

\smallskip

\begin{theorem}\label{submain}
Let $R$ be a local Gorenstein ring with infinite residue field,
let $I$ be an $R$-ideal with $g={\rm ht} \, I$ and $\ell=
\ell(I)$, and let $J$ be a minimal reduction of $I$. Assume $I$
satisfies $G_{\ell}$ and ${\rm depth} \, R/I^j \geq {\rm dim} \,
R/I -j+1 \, $ for $1 \leq j \leq \ell-g$. Then
\[
{\rm core}(I) \subset J^{n+1} : \displaystyle \sum_{y \in I} (J,
y)^n
\]
for every $n \geq 0$.
\end{theorem}
\begin{proof}
We may assume that $\ell > 0$ and $n > 0$. We use the
notation of Lemma~\ref{lifting} with $t=n+1$ and $H$ the
intersection of all primary components of $J^{n+1}$ of height $<
\ell$. Notice that the assumptions of Lemma~\ref{lifting} are
satisfied by Remark~\ref{lineartype}. Write $L=\displaystyle
\sum_{y \in I} (J, y)^n$.

We first prove
\begin{equation}\label{eq}
{\rm core}(I) \subset H \colon L,
\end{equation}
or equivalently, $({\rm core}(I))_{\mathfrak p} \subset (H \colon
L)_{\mathfrak p}$ for every prime ideal ${\mathfrak p}$ with ${\rm
dim}\, R_{\mathfrak p} < \ell$. Indeed by Remark~\ref{lineartype},
$I_{\mathfrak p}=J_{\mathfrak p}$  and hence $L_{\mathfrak
p}=J^n_{\mathfrak p}$. Thus $({\rm core}(I))_{\mathfrak p} \subset
J_{\mathfrak p} \subset J^{n+1}_{\mathfrak p} \colon
J^n_{\mathfrak p} = H_{\mathfrak p} \colon L_{\mathfrak p}$, which
shows $(\ref{eq})$.

Next let ${\mathfrak a} \in {\mathcal A}$. We prove that
\begin{equation}\label{eq2}
{\rm core}(I) \subset (J^{n+1}, {\mathfrak a}) \, \colon L.
\end{equation}
Let `${}^{\tratto}$' denote images in $\overline{R} =
R/{{\mathfrak a} \, \colon I}$. According to
Remark~\ref{lineartype}, ${\rm ht} \, J \colon I \geq \ell$ and
hence ${\mathfrak a} \, \colon I$ is a geometric
$(\ell-1)$-residual intersection of $I$. Thus by \cite[2.9$({\it
a})$ and 1.7$({\it a,c})$]{U}, $({\mathfrak a} \, \colon I) \cap I
={\mathfrak a}$, ${\rm ht}\, \overline{I} \geq 1$ and
$\overline{R}$ is Cohen-Macaulay. Furthermore
$\ell(\overline{I})\leq 1$. Hence $\ell(\overline{I}) = {\rm ht}\,
\overline{I}=1$ and $\overline{J}$ is a minimal reduction of
$\overline{I}$. As $\overline{J}$ is generated by a single regular
element, $\overline{J}^{j+1} \colon \displaystyle
\sum_{\overline{y} \in \overline{I}} (\overline{J},
\overline{y})^j$ form a decreasing sequence of ideals. Thus
Theorem~\ref{heightone}(${\it a}$) implies that ${\rm
core}(\overline{I}) \subset \overline{J}^{n+1} \colon
\overline{L}$. On the other hand by \cite[4.5]{CPU}, ${\rm
core}(\overline{I}) = (\overline{\alpha}_1) \cap \ldots \cap
(\overline{\alpha}_t)$ for some integer $t$ and $t$ general
principal ideals $(\alpha_1), \ldots, (\alpha_t)$ in $I$. Notice
that $({\mathfrak a}, \alpha_i)$ are reductions of $I$, hence
${\rm core}(I) \subset \displaystyle\bigcap_{i=1}^t ({\mathfrak
a}, \alpha_i)$. Therefore $ \, \overline{{\rm core}(I)} \subset
\displaystyle\overline{\bigcap_{i=1}^t({\mathfrak a}, \alpha_i)}
\subset \bigcap_{i=1}^t (\overline{\alpha}_i) = {\rm
core}(\overline{I})$. As  ${\rm core}(\overline{I}) \subset
\overline{J}^{n+1} \colon \overline{L}$ we obtain
\begin{eqnarray*}
{\rm core}(I)  & \subset & (J^{n+1}, {\mathfrak a}\,  \colon I) \,\colon L \\
& = & (J^{n+1}, ({\mathfrak a} \, \colon I) \cap I)\, \colon L \\
& = & (J^{n+1}, {\mathfrak a})\, \colon L,
\end{eqnarray*}
\noindent which proves $(\ref{eq2})$.

Finally, combining $(\ref{eq})$, $(\ref{eq2})$ and
Lemma~\ref{lifting} we deduce
\begin{eqnarray*}
 {\rm core}(I) & \subset &
 (H \cap \bigcap _{{\mathfrak a} \in {\mathcal A}}
 (J^{n+1}, {\mathfrak a})) \, \colon L  \\
 & \subset & J^{n+1} \colon L = J^{n+1} \colon \displaystyle \sum_{y \in I} (J,
 y)^n,
\end{eqnarray*}
as claimed. \QED
\end{proof}

\smallskip

We are now ready to assemble the proof of the main theorem.

\begin{theorem}\label{main}
Let $R$ be a local Gorenstein ring with infinite residue field
$k$, let $I$ be an $R$-ideal with $g={\rm ht} \, I > 0$ and $\ell=
\ell(I)$, and let $J$ be a minimal reduction of $I$ with
$r=r_J(I)$. Assume $I$ satisfies $G_{\ell}$, ${\rm depth} \, R/I^j
\geq {\rm dim} \, R/I -j+1$ for $1 \leq j \leq \ell-g$, and either
${\rm char} \,  k =0 $ or ${\rm char} \, k  > r -\ell + g $. Then
\[
{\rm core}(I)= J^{n+1} : I^n
\]
for every $n \geq {\rm max} \{ r - \ell + g, 0 \}$.
\end{theorem}
\begin{proof}
To prove the containment $J^{n+1} : I^n \subset {\rm
core}(I)$ we show that if $K$ is an arbitrary minimal reduction of
$I$ then $J^{n+1} : I^n \subset K$. By Corollary~\ref{balanced},
$J^{n+1} : I^n = K^{j+1} : I^j$ for $j \gg 0$. However, $K^{j+1} :
I^j \subset K^{j+1} : K^j = K$, since $g > 0$ and ${\rm
gr}_{K}(R)$ is Cohen-Macaulay by Remark~\ref{lineartype}. Thus
indeed $J^{n+1} \colon I^n \subset K$.

To show the inclusion ${\rm core}(I) \subset J^{n+1} \colon I^n$,
we may assume $n = {\rm max}\{ r - \ell +g , 0 \}$ again by
Corollary~\ref{balanced}. But then $I^n= \displaystyle \sum_{y \in
I} (J, y)^n $ according to our assumption on the characteristic,
and the asserted containment follows from
Theorem~\ref{submain}. \QED
\end{proof}

\smallskip

In the next corollary we prove that the core of $I$ is integrally
closed whenever
%${\rm Proj}(R[It])$
$R[It]$ is regular in codimension one. Similar results can be
found in \cite[3.12]{HS}, \cite[2.11]{CPU2}, \cite[5.5.3]{HyS},
where essentially the normality of the ideal and the
Cohen-Macaulayness of the Rees ring were required.

%\begin{Corollary}
%If in addition to the assumptions of Theorem~\ref{main}, $R$ is
%normal and analytically unramified
%and $R[It]$ satisfies Serre's condition $R_1$,
%then ${\rm core}(I)$ is integrally closed.
%\end{Corollary}
%\demo Let $J$ be a minimal reduction of $I$, set $A=R[Jt,t^{-1}]
%\subset B= R[It,t^{-1}]$, write $\overline{B} \, $ for the
%integral closure of $B$ in $R[t,t^{-1}]$, and $L= \oplus_i
%(J^{n+i+1} \colon I^n) t^i$ where $n \geq {\rm max} \{ r - \ell +
%g, 0 \}$. Notice that $L \subset A \subset \overline{B} \, $ by
%the proof of Proposition~\ref{canonical}, and $[L]_0 ={\rm
%core}(I)$ by Theorem~\ref{main}. According to the assumptions of
%the corollary, $\overline{B} \, $ is a normal domain and the
%$A$-module $\overline{B}/B$ has grade $\geq 2$. Thus
%$\omega_{\overline{B}} \cong \omega_B$, whereas $\omega_B(g) \cong
%L$ by Remark~\ref{remarkcanonical}(2). Therefore $L$ is a
%$\overline{B}$-module satisfying $S_2$, and hence is a divisorial
%ideal of the normal domain $\overline{B}$. Thus $L$ is an
%integrally closed $\overline{B}$-ideal, which forces $[L]_0$ to be
%integrally closed $R$-ideal.
%%as an ideal in $[\overline{B}]_0=R$.
% \QED

\begin{corollary}
If in addition to the assumptions of Theorem~\ref{main},
$R[It]$ satisfies Serre's condition $R_1$,
then ${\rm core}(I)$ is integrally closed.
\end{corollary}
\begin{proof}
Let $A=R[Jt,t^{-1}] \subset B= R[It,t^{-1}]$, and for  $n \geq
{\rm max} \{ r - \ell + g, 0 \}$ write $L= \oplus_i (J^{n+i+1}
\colon I^n) t^i$. Notice that $L \subset A \subset B \, $ by the
proof of Proposition~\ref{canonical}, $[L]_0 ={\rm core}(I)$ by
Theorem~\ref{main}, and $\omega_B(g) \cong L$ by
Remark~\ref{remarkcanonical}(2). Therefore $L$ is a $B$-module
satisfying $S_2$, and hence an unmixed $B$-ideal of height $1$. As
$B$ satisfies $R_1$, $L$ is an integrally closed $B$-ideal, which
forces $[L]_0$ to be integrally closed as an $R$-ideal.
%as an ideal in $[\overline{B}]_0=R$.
\QED
\end{proof}

\smallskip

The next result has been shown in \cite[5.3.1]{HyS} for equimultiple
ideals.

\begin{remark}{\rm
If in addition to the assumptions of Theorem~\ref{main}, $R$ is a
regular local ring essentially of finite type over a field of
characteristic zero and $R[It]$ has only rational singularities,
then
\[
{\rm core}(I)= {\rm adj}(I^{g})= I \, {\rm adj}(I^{g-1}),
\]
where ${\rm adj}$ denotes adjoint ideals in the sense of
\cite[1.1]{L}.

Indeed, notice that ${\rm core}(I) = [\omega_B]_g t^{-g}= I
[\omega_B]_{g-1} t^{-g+1}$ by Theorem~\ref{main} and
Remarks~\ref{remarkcanonical}$($2$)$ and~\ref{Jout}, whereas $
[\omega_B]_i t^{-i}= {\rm adj}(I^i)$ for every $i \geq 0$
according to \cite[the proof of 3.5]{Hy}.}
\end{remark}

\smallskip

\begin{remark}{\rm
Even without the assumption on the characteristic of the residue
field in Theorem~\ref{main}, Theorem~\ref{submain} and
the proof of Theorem~\ref{main} still show that
\[
J^{n+1} \colon I^n \subset {\rm core} (I) \subset J^{n+1} \colon
\sum_{y \in I} (J, y)^n
\]
for $n \geq {\rm max} \{ r - \ell + g , 0 \}$. In particular
if $ \mu(I) \leq \ell +1$ then
\[
{\rm core} (I) = J^{n+1} \colon I^n.
\]
}
\end{remark}

\bigskip

In general however, the formula of Theorem~\ref{main} is no longer
valid if $\mu(I) > \ell +1$ and $0 < {\rm char} \, k \leq r - \ell
+ g $. We illustrate this with a class of examples in which the
ambient ring is a domain:

\begin{example}\label{charp}
{\rm Let $k$ be an infinite field of characteristic $p > 0$, let
$q > p$ be an integer not divisible by $p$, consider the numerical
semigroup ring $R=k[\![ t^{p^2}, t^{pq}, t^{pq +q} ]\!] \subset
k[\![ t ]\!]$, and let $I = {\mathfrak m}$ be the maximal ideal of
$R$. Now $R$ is a one-dimensional local Gorenstein domain, and one
has the {\it proper} containment ${\rm core} (I) \supsetneq
J^{n+1} \colon I^n$ for any minimal reduction $J$ of $I$ and any
$n \geq r(I)$. In fact ${\rm core}(I)= (t^{p^3}, {\mathfrak
m}^{2p-1})$, whereas $J^{n+1} \colon I^n= {\mathfrak m}^{2p-1}$.

To prove these claims consider the presentation $R \cong k[\![ X,
Y, Z]\!] / (Y^p -X^q, Z^p -X^qY) $ where $X, Y, Z$ are mapped to
$x=t^{p^2}, y=t^{pq}, z=t^{pq+q}$, respectively. Clearly $R$ is
Gorenstein. Comparing multiplicities one sees that ${\rm
gr}_{\mathfrak m}(R) \cong k[X,Y,Z]/(Y^p,Z^p)$. Thus $r({\mathfrak
m})=2(p-1)$. Furthermore the leading form $x^{*}$ of $x$ in ${\rm
gr}_{\mathfrak m}(R)$ is a regular element, in particular $x$
generates a minimal reduction of ${\mathfrak m}$.

We apply Theorem~\ref{heightone}(${\it a}$) taking $(x)$ as the
minimal reduction of $I$ and using the definition of $y_i$ and $s$
as in that theorem. Since $(y^p,z^p)=(x^q,x^qy) \subset (x^p)$ and
${\rm char} \, k =p$, one has $y_i^p \in (x^p)$ for every $i$, and
therefore $s \leq p-1 < r({\mathfrak m})$. As $p-1 < {\rm char} \,
k$, Theorem~\ref{heightone}(${\it a}$) implies
\[
{\rm core}({\mathfrak m}) = (x^p) \colon {\mathfrak m}^{p-1} =
(x^{n+1}) \colon x^{n+1-p}{\mathfrak m}^{p-1}.
\]
On the other hand, $J^{n+1} \colon {\mathfrak m}^n = (x^{n+1})
\colon {\mathfrak m}^n$ by Lemma~\ref{help}$({\it a})$, and
\[
(x^{n+1}) \colon {\mathfrak m}^n \subsetneq (x^{n+1}) \colon
x^{n+1-p}{\mathfrak m}^{p-1}
\]
because $R$ is Gorenstein and $n \geq r({\mathfrak m}) > p-1$.
Thus ${\rm core}({\mathfrak m}) \varsupsetneq J^{n+1} \colon
{\mathfrak m}^n$.

Now, to compute these ideals let us write `${}^{\tratto}$' for
images in $\overline{R} = R/(x^p)$. One has ${\rm
gr}_{\overline{\mathfrak m}}(\overline{R}) \cong {\rm
gr}_{\mathfrak m}(R)/ (x^{*p}) \cong k[X,Y,Z]/(X^p,Y^p,Z^p)$ by
the regularity of $x^*$, hence
\[
{\rm core}({\mathfrak m}) = (x^p) \colon {\mathfrak m}^{p-1} =
(x^p, {\mathfrak m}^{2p-1})=(t^{p^3}, {\mathfrak m}^{2p-1}).
\]
Likewise one sees that
\[
J^{n+1} \colon {\mathfrak m}^n = (x^{n+1}) \colon {\mathfrak m}^n
= (x^{n+1}, {\mathfrak m}^{2p-1}) ={\mathfrak m}^{2p-1}.
\]}
%where the last equality holds because $n \geq r({\mathfrak m}) =
%2p -2$.}
\end{example}

\bigskip

In the next example we show that the $G_{\ell}$ condition cannot
be removed from Theorem~\ref{main}.

\begin{example}
{\rm Let $k$ be an infinite field, write $R = k[\![X, Y, Z,
W]\!]/(X^2+Y^2+Z^2, ZW)$, let $x, y, z, w$ denote the images of
$X, Y, Z, W$ in $R$, and consider the $R$-ideal $I=(x, y, z)$.
Notice that $R$ is a local Gorenstein ring, ${\rm ht}\, I =1$,
$\ell(I) = 2$, $R/I$ is Cohen-Macaulay, but $I$ does not satisfy
$G_2$. Let $J=(x,y)$. The ideal $J$ is a minimal reduction of $I$
with $r_J(I)=1$. One has ${\rm core}(I)=I^2 \varsubsetneq J^2
\colon I$. The same holds if one replaces $J$ by a general minimal
reduction of $I$.

Indeed, the special fiber ring ${\rm gr}_I(R) \otimes_R k$ is
defined by a single quadric. Hence $\ell(I) = 2$, and $r_K(I) = 1$
for every minimal reduction $K$ of $I$, which gives $I^2 \subset
{\rm core}(I)$. On the other hand,  $(x,y)$, $(x,z)$ and $(y,z)$
are minimal reductions of $I$, thus ${\rm core}(I) \subset (x,y)
\cap (x,z) \cap (y,z) = I^2$. Therefore ${\rm core}(I) = I^2$. To
conclude notice that $I^2 \varsubsetneq  (I^2, xw, yw)=J^2 \colon
I$.}
\end{example}

\bigskip

Finally, the formula of Theorem~\ref{main} does not hold for $g=0$
even if $\ell > 0$:

\begin{example}
{\rm Let $k$ be an infinite field, let $\Delta_i \in
k[\![X,Y,Z]\!]$ be the maximal minor of the matrix
\[
\left(
\begin{array}{cccc}
X & Y & 0 & Z \\
Y & 0 & Z & X \\
0 & Z & X & Y
\end{array} \right)
\]
obtained by deleting the $i^{\rm th}$ column, set $R=
k[\![X,Y,Z]\!]/(\Delta_1, \Delta_2)$ and define $J=\Delta_3R$,
$I=(\Delta_3, \Delta_4)R$.  Then $R$ is a local Gorenstein ring,
${\rm ht} \, I =0 $, $\ell(I)=1$, $I$ satisfies $G_1$, $R/I$ is
Cohen-Macaulay, and $J$ is a minimal reduction of $I$ with
$r=r_J(I)=2$. However,  ${\rm core}(I) \subsetneq J^{n+1} \colon
I^n$ for every $n \geq 1={\rm max} \{ r - \ell +g, 0\}$.

Indeed, \cite[5.1]{U2} and \cite[3.6]{PU} show that $J$ is a
minimal reduction of $I$ with $r=2$. Writing ${\mathfrak
m}=(X,Y,Z)R$ one has $J:I={\mathfrak m}$. As $r=2= \ell -g +1$,
\cite[2.6(2)]{CPU2} then gives ${\rm core}(I)= {\mathfrak
m}J={\mathfrak m}I$. On the other hand, a computation shows that
${\mathfrak m}J \subsetneq J^2 \colon I = J^3 \colon I^2$. The
assertion now follows since $J^{n+1} \colon I^n$ form an
increasing sequence of ideals for $n \geq r=2$. }
\end{example}

\begin{remark}\label{ConjJout}
{\rm If in addition to the assumptions of Theorem~\ref{main},
${\rm gr}_I(R)$ is Cohen-Macaulay, then
\[
{\rm core}(I)= J (J^n :I^n) = I (J^n : I^n)
\]
for every $n \geq {\rm max} \{ r - \ell + g, 0 \}$.

This assertion follows from Theorem~\ref{main} and
Remark~\ref{Jout}.}
\end{remark}

%We conjecture that Remark~\ref{ConjJout} holds without the
%Cohen-Macaulay assumption on the associated graded ring:

%\smallskip

%\begin{Conjecture}
%In the setting of Theorem~\ref{main},
%\[
%{\rm core}(I) = J(J^n \colon I^n) = I(J^n \colon I^n).
%\]
%\end{Conjecture}

\medskip

In the light of Theorem~\ref{main} and
Remark~\ref{remarkcanonical}$($2$)$ the formula of
Remark~\ref{ConjJout} is equivalent to saying that the canonical
module of $B=R[It,t^{-1}]$ as a graded module over
$A=R[Jt,t^{-1}]$ has no homogeneous minimal generators in degree
$g$. In fact one could ask whether the $A$-module $\omega_B$ is
generated in degrees $ \leq g-1$, or equivalently, whether
\[
J^{i+n} \colon I^n = J^i(J^n \colon I^n) = I^i(J^{n} \colon I^n)
\]
for every $i \geq 0$. This stronger statement still holds if ${\rm
gr}_I(R)$ is Cohen-Macaulay, according to Remark~\ref{Jout}. For
$I$ an equimultiple ideal it also holds under the weaker condition
that the $S_2$-ification of $B$ is Cohen-Macaulay (i.e.,
$\omega_B$ is Cohen-Macaulay), or under the assumption that $R$ is
a domain essentially of finite type over a field of characteristic
zero and ${\rm Proj}(R[It])$ is smooth. For the latter case one
uses a Kodaira type vanishing theorem due to
Cutkosky (\cite[A2]{L}).

\bigskip

We finish with an example showing that the equality of
Remark~\ref{ConjJout} does not hold in general, even for
${\mathfrak m}$-primary ideals of reduction number two in
two-dimensional regular local rings. This example was constructed
by Heinzer, Johnston and Lantz for a slightly different purpose
(\cite[5.4]{HJL}).

\begin{example}
{\rm Let $k$ be an infinite field with ${\rm char}\, k \neq 2$,
let $ R=k[\![X,Y]\!]$, and $I=(X^7, Y^7, X^3Y^5+X^5Y^3, X^4Y^6)$.
One has $r=r_J(I)=2$, but $ J (J^n :I^n) \subset I (J^n :I^n)
\subsetneq J^{n+1} : I^n={\rm core}(I)$  for every minimal
reduction $J$ of $I$ and every $n \geq 2={\rm max} \{ r - \ell +g,
0\}$.

Indeed, a computation shows that ${\rm depth}\,{\rm gr}_I(R)=1$ and
$r(I)=2$, see \cite[5.4]{HJL}, and then
$r_J(I)=2$ for every minimal reduction $J$ of $I$ according to
\cite[2.1]{Huc} or \cite[1.2]{Trung}. Now
%Theorem~\ref{main} and Remark~\ref{remarkcanonical}$($2$)$ shows
%that $ I (J^n :I^n)
%\sebset  J^{n+1} : I^n={\rm core}(I)$  and that these ideals are
%independent of $J$ and $n\geq2$.
Remark~\ref{remarkcanonical}$($2$)$ shows that $I (J^n :I^n)$ is
independent of $J$ and $n \geq 2$, and Theorem~\ref{main} gives
$J^{n+1} : I^n={\rm core}(I)$. Finally, taking $J=(X^7,Y^7)$ one
easily computes that $I (J^2 :I^2) \subsetneq J^3 :I^2$.}
\end{example}


\begin{thebibliography}{99}

\bibitem{AH2}{I.M. Aberbach and C. Huneke, A theorem of
Brian\c{c}on-Skoda type for regular local rings containing a
field, Proc. Amer. Math. Soc. {\bf 124} (1996), 707-713.}


\bibitem{CPU}{A. Corso, C. Polini and B. Ulrich, The structure of
the core of ideals, Math. Ann. {\bf 321} (2001), 89-105.}

\bibitem{CPU2}{A. Corso, C. Polini and B. Ulrich, Core and residual
intersections of ideals, Trans. Amer. Math. Soc. {\bf 354} (2002),
2579-2594.}

\bibitem{CPU3}{A. Corso, C. Polini and B. Ulrich,
The core of projective dimension one modules,
Manuscripta Math. {\bf 111} (2003), 427-433.}

\bibitem{HJL}{W. Heinzer, B. Johnston and D. Lantz, First
coefficient domains and ideals of reduction number one, Comm.
Algebra {\bf 21} (1993), 3797-3827.}

\bibitem{HSV}{J. Herzog, A. Simis and W.V. Vasconcelos, Koszul
homology and blowing-up rings, in {\it Commutative Algebra,
Proceedings: Trento 1981} (Greco/Valla eds.), Lecture Notes in
Pure and Applied Mathematics {\bf 84}, Marcel Dekker, New York,
1983, 79-169.}

\bibitem{HVV}{J. Herzog, W.V. Vasconcelos and R. Villarreal,
Ideals with sliding depth, Nagoya Math. J. {\bf 99} (1985),
159-172.}

\bibitem{HH}{M. Hochster, C. Huneke,
Indecomposable canonical modules and connectedness, Contemp. Math.
{\bf 159} (1994), 197-208.}

\bibitem{Huc}{S. Huckaba,
Reduction numbers for ideals of higher analytic spread, Math.
Proc. Camb. Phil. Soc. {\bf 102} (1987), 49-57.}

\bibitem{H1}{C. Huneke, Linkage and Koszul homology of ideals,
Amer. J. Math. {\bf 104} (1982), 1043-1062.}

\bibitem{HS}{C. Huneke and I. Swanson, Cores of ideals in $2$-dimensional
regular local rings, Michigan Math. J. {\bf 42} (1995), 193-208.}

\bibitem{HT}{C. Huneke and N. V. Trung, On the core of ideals,
preprint.}

\bibitem{Hy}{E. Hyry, Coefficient ideals and the Cohen-Macaulay property of
Rees algebras, Proc. Amer. Math. Soc. {\bf 129} (2001),
1299-1308.}

\bibitem{HyS}{E. Hyry and K. Smith, On a non-vanishing conjecture
of Kawamata and the core of an ideal, Amer. J. Math. {\bf 125}
(2003), 1349-1410.}

\bibitem{Lip}{J. Lipman, Stable ideals and Arf rings, Amer. J. Math.
{\bf 93} (1971), 649-685.}

\bibitem{L}{J. Lipman, Adjoints of ideals in regular local rings,
Math. Res. Lett. {\bf 1} (1994), 739-755.}


\bibitem{PU}{C. Polini and B. Ulrich, Linkage and reduction numbers,
Math. Ann. {\bf 310} (1998), 631-651.}


\bibitem{RS}{D. Rees and J.D. Sally, General elements and joint
reductions, Michigan Math. J. {\bf 35} (1988), 241-254.}

\bibitem{Trung}{N.V. Trung, Reduction exponent and degree
bound for the defining equations of graded rings, Proc. Amer.
Math. Soc. {\bf 101} (1987), 229-236.}

\bibitem{U}{B. Ulrich, Artin-Nagata properties and reductions of ideals,
Contemp. Math. {\bf 159} (1994), 373-400.}

\bibitem{U2}{B. Ulrich, Ideals having the expected reduction number,
Amer. J. Math. {\bf 118} (1996), 17-38.}

\end{thebibliography}
\end{document}